\documentclass{amsart}
\usepackage[centertags]{amsmath}
\usepackage{amsfonts}
\usepackage{amssymb}
\usepackage{amsthm}
\usepackage[ansinew]{inputenc}
\usepackage[all]{xy}
\usepackage{graphicx}
\newtheorem{thm}{Theorem}[section]
\newtheorem{cor}[thm]{Corollary}
\newtheorem{lem}[thm]{Lemma}
\newtheorem{prop}[thm]{Proposition}
\theoremstyle{definition}
\newtheorem{defn}[thm]{Definition}
\theoremstyle{remark}
\newtheorem{rem}[thm]{Remark}
\newtheorem{exm}[thm]{Example}
\numberwithin{equation}{section}

\newcommand{\set}[1]{\left\{#1\right\}}

\newcommand{\To}{\longrightarrow}

\newcommand{\C}[1]{\mathbf{#1}} 

\def\r{\rightarrow} 
\def\l{\leftarrow} 

\def\hom{\operatorname{Hom}}
\def\ext{\operatorname{Ext}}

\def\F{\mathcal{F}} 

\def\st{\stackrel} 

\newcommand{\Ab}{\mathrm{\mathbf{Ab}}}

\def\Z{\mathbb{Z}}

\numberwithin{equation}{section}

\begin{document}

\title{On the functoriality of cohomology of categories}%
\author{Fernando Muro}%
\address{Departamento de Geometría y Topología de la Universidad de Sevilla}%
\email{fmuro@us.es}%

\thanks{The author was partially supported
by the MCyT grant BFM2001-3195-C02-01 and the MECD FPU fellowship
AP2000-3330.}%

\begin{abstract}
In this paper we show that the Baues-Wirsching complex used to
define cohomology of categories is a $2$-functor from a certain
$2$-category of natural systems of abelian groups to the
$2$-category of chain complexes, chain homomorphism and relative
homotopy classes of chain homotopies. As a consequence we derive
(co)localization theorems for this cohomology.
\end{abstract}
\maketitle

\section{Introduction}

Baues-Wirsching cohomology of a small category $\C{C}$ with
coefficients in a natural system $D$ on $\C{C}$ was defined in
\cite{csc} as the cohomology of a certain cochain complex
$F^*(\C{C},D)$. This cohomology generalizes some other
cohomologies previously known, as for example
\begin{itemize}
\item the Hochschild-Mitchell cohomology of $\C{C}$ with
coefficients in a functor
$D\colon\C{C}^\mathrm{op}\times\C{C}\r\C{Ab}$ (\cite{rso}),

\smallskip

\item the cohomology of the classifying space $B\C{C}$ with
local coefficients $D$,

\smallskip

\item the Mac Lane cohomology of a ring (\cite{ham}, \cite{cat}).
\end{itemize}

The Baues-Wirsching complex $F^*(\C{C},D)$ as well as its
cohomology $H^*(\C{C},D)$ are known to be functors on a certain
category $\C{Nat}$ of pairs $(\C{C},D)$. It is also known that
equivalences of categories induce homotopy equivalences in
Baues-Wirsching complexes and isomorphisms in cohomology groups,
however this does not follow immediately from the fact that $F^*$
and $H^*$ are functors in $\C{Nat}$. More precisely, it is not
known the behaviour of $F^*(\C{C},D)$ and $H^*(\C{C},D)$ with
respect to natural transformations between functors in the first
variable. The goal of this paper is to shed some light on that
issue.

We define a new $2$-category $\C{Nat}_\F$ of pairs $(\C{C},D)$
containing $\C{Nat}$ and prove that $F^*$ is in fact a $2$-functor
in $\C{Nat}_\F$ (Theorem \ref{2F}). The $2$-morphisms in the
category of cochain complexes will be homotopy classes of
homotopies relative to the boundary of the cylinder. The category
obtained from $\C{Nat}_\F$ by taking sets of connected components
on morphism categories turns out to be a quotient of $\C{Nat}$
(Proposition \ref{full}) and the cohomology functor $H^*$ factors
through this quotient category (Corollary \ref{elcor}).

As an application of these results we obtain localization and
colocalization theorems for Baues-Wirsching cohomology (Theorems
\ref{locthm} and \ref{colocthm}). We also give some examples of
how these (co)localization theorems can be used to carry out
computations in cohomology of categories.

Pirashvili and Waldhausen defined in \cite{mlhthh} the homology of
a small category $\C{C}$ with coefficients in a functor
$D\colon\C{C}^\mathrm{op}\times\C{C}\r\C{Ab}$ by using a complex
$F_*(\C{C},D)$ which is similar, and in some sense dual, to the
Baues-Wirsching complex. They also proved that this homology
extends Mac Lane homology of rings (\cite{ham}), which is
isomorphic to topological Hochschild homology in the sense of
\cite{thh}. We claim that the definition of the
Pirashvili-Waldhausen homology $H_*(\C{C},D)$ can be extended to
natural systems $D$ as coefficient objects and the functorial
properties of $F_*$ and $H_*$ are similar to those described here
for $F^*$ and $H^*$, in particular (co)localization theorems
should hold.

\section{Notation and conventions}

Compelled by the necessarily intricate notation of this paper, we
have decided to include at the beginning a section to fix the
meaning of some symbols. Additional notation will also appear in
the development of the paper, but it will not contradict in any
case that introduced here.

\bigskip

\begin{center}
\begin{tabular}{|c|l|}
  \hline \textit{symbols} & \textit{meaning}\\
  \hline \hline $\C{C},\C{D},\C{E}$ & small categories. \\
  \hline $X,Y$ & objects in those categories. \\
  \hline $f,g,h,k,\sigma$ & morphisms in those categories. \\
  \hline $\varphi,\psi,\xi,\zeta$ & functors between those categories. \\
  \hline $\alpha, \beta, \gamma, \varepsilon, t, s$ & natural transformations. \\
  \hline $\r$ & arrow for morphisms in ordinary categories, \\ &$1$-morphisms in $2$-categories, functors, and $2$-functors. \\
  \hline $\Rightarrow$ & arrow for $2$-morphisms and natural transformations. \\
  \hline $1$ & identity morphism or identity $2$-morphism, \\ & a subscript will clarify the meaning in ambiguous cases.\\
  \hline $\Ab$ & the category of \\ & \hspace{10pt}\textit{objects:} abelian groups, \\ & \hspace{10pt}\textit{morphisms:} homomorphisms.\\
  \hline $\C{Cat}$ & the category of \\ & \hspace{10pt}\textit{objects:} small categories, \\ & \hspace{10pt}\textit{morphisms:} functors.\\
  \hline $\C{Cat}_2$ & the standard $2$-category of \\ &\hspace{10pt}\textit{objects:} small categories, \\ & \hspace{10pt}\textit{$1$-morphisms:} functors, \\ & \hspace{10pt}\textit{$2$-morphisms:} natural transformations.\\
  \hline $D, E, G$ & natural systems, see Definition \ref{cc}. \\
  \hline $c$ & cochain in cohomology of categories.\\
  \hline $A^*, B^*, C^*$ & cochain complexes of abelian groups. \\
  \hline $d$ & the differential in all cochain complexes. \\
  \hline $p,q,h,r$ & graded morphisms between graded abelian
  groups.\\ \hline
\end{tabular}
\end{center}

\bigskip

These symbols can be altered by adding superscripts or subscripts.

In $2$-categories the word ``morphism'' will be synonym of
``$1$-morphism''. All categories can be regarded as $2$-categories
with only the trivial $2$-morphisms.

The symbol $\bullet$ stands for an unspecified object in an
arbitrary category. It can appear several times in a diagram,
however in general it will stand for a different object each time.

The composite of morphisms, say
$\bullet\st{f}\r\bullet\st{g}\r\bullet$, or functors will be
indicated by juxtaposition $gf$, as well as the vertical
composition of $2$-morphisms or natural transformations
$\beta\alpha\colon\varphi\Rightarrow\xi$ as in the diagram
$$
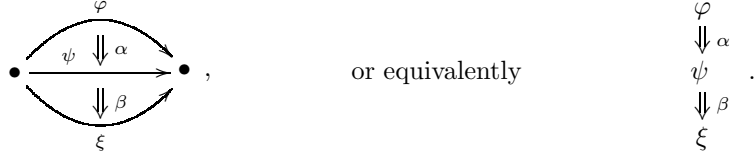
\xymatrix@C=50pt@R=10pt{&&&\varphi\ar@{=>}[d]^{\;\alpha}\\\bullet\ar@/^20pt/[r]^{\varphi}_{\;}="a"\ar[r]^<(.25){\psi}^{}="b"_{\;}="c"\ar@/_20pt/[r]_{\xi}^{}="d"
&**[r]\bullet\text{\hspace{5pt},}\ar@{=>}^{\;\alpha}"a";"b"\ar@{=>}^{\;\beta}"c";"d"&\text{or
equivalently}&**[r]\psi\text{\hspace{15pt}.}\ar@{=>}[d]^{\;\beta}\\&&&\xi}$$

We will use the symbol $*$ for the horizontal composition of
$2$-morphisms or natural transformations
$\beta*\alpha\colon\xi\varphi\Rightarrow\zeta\psi$ as in the
following diagram
$$
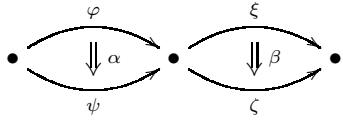
\xymatrix@C=50pt{\bullet\ar@/^12pt/[r]^\varphi_{\;}="a"\ar@/_12pt/_\psi[r]^{\;}="b" &\bullet\ar@{=>}^{\;\alpha}"a";"b"
\ar@/^12pt/[r]^\xi_{\;}="c"\ar@/_12pt/[r]_\zeta^{\;}="d"
&\bullet\ar@{=>}^{\;\beta}"c";"d"}$$

\section{Factorization categories}

Factorization categories are the source of the coefficient objects
for the Baues-Wirsching cohomology of small categories, see
Definition \ref{cc}. A thorough study of their properties is
essential to study in depth the functoriality of the
Baues-Wirsching complex in Section 5.

\begin{defn}\label{fc}
The \emph{factorization category} $\F\C{C}$ of a small category
$\C{C}$ has
\begin{itemize}
\item[] \textit{objects:} morphisms in $\C{C}$,

\smallskip

\item[] \textit{morphisms:} $(h,k)\colon f\r g$ are pairs of
morphisms in $\C{C}$ such that $kfh=g$, that is commutative
diagrams in $\C{C}$
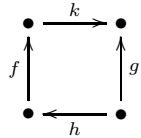
$$\xymatrix{\bullet\ar[r]^k&\bullet\\\bullet\ar[u]^f&\bullet\ar[u]_g\ar[l]^h}$$
and composition is defined by $(h',k')(h,k)=(hh',k'k)$.
\end{itemize}
\end{defn}

One can easily check that factorization categories define a
functor
\begin{equation*}
\F\colon\C{Cat}\To\C{Cat}.
\end{equation*}
This functor is defined on morphisms as follows: a functor
$\varphi\colon\C{C}\r\C{D}$ is sent to another one
\begin{equation*}
\F(\varphi)\colon\F\C{C}\To\F\C{D},
\end{equation*}
which is given
\begin{itemize}
\item[] \textit{on objects:} by $\F(\varphi)(f)=\varphi(f)$,

\smallskip

\item[] \textit{on morphisms:} by $\F(\varphi)(h,k)=(\varphi(h),\varphi(k))$.
\end{itemize}

Notice that the functor $\F$ preserves products, therefore we can
consider the $2$-category $\C{Cat}_\F$ obtained from $\C{Cat}_2$
by applying the functor $\F$ to morphism categories. Let us make
explicit the structure of $\C{Cat}_\F$:
\begin{itemize}
\item[(3.A)\hspace{8pt}] \textit{objects:} are small categories;

\smallskip

\item[] \textit{$1$-morphisms:} $\alpha\colon\C{C}\r\C{D}$ are
actually natural transformations
$\alpha\colon\varphi\Rightarrow\psi$ between functors
$\varphi,\psi\colon\C{C}\r\C{D}$, and composition $\beta\alpha$ in
$\C{Cat}_\F$ is horizontal composition $\beta*\alpha$ of natural
transformations;

\smallskip

\item[] \textit{$2$-morphisms:}
$(\varepsilon,\gamma)\colon\alpha\Rightarrow\beta$ are natural
transformations such that $\gamma\alpha\varepsilon=\beta$, that is
commutative diagrams of natural transformations
$$\xymatrix{\psi\ar@{=>}[r]^\gamma&\zeta\\\varphi\ar@{=>}[u]^\alpha&\xi\ar@{=>}[u]_\beta\ar@{=>}[l]^\varepsilon}$$
vertical composition of $2$-morphisms is given by
$(\varepsilon',\gamma')(\varepsilon,\gamma)=(\varepsilon\varepsilon',\gamma'\gamma)$,
and the horizontal composition of $(\varepsilon,\gamma)$ and
$(\varepsilon',\gamma')$ as in the following diagram in
$\C{Cat}_\F$
$$\xymatrix@C=70pt{\C{C}\ar@/^15pt/[r]^{\alpha}_{\;}="a"\ar@/_15pt/_{\beta}[r]^{\;}="b" &\C{D}\ar@{=>}^{\;(\varepsilon,\gamma)}"a";"b"
\ar@/^15pt/[r]^{\alpha'}_{\;}="c"\ar@/_15pt/[r]_{\beta'}^{\;}="d"
&\C{E}\ar@{=>}^{\;(\varepsilon',\gamma')}"c";"d"}$$ is
$(\varepsilon',\gamma')*(\varepsilon,\gamma)=(\varepsilon'*\varepsilon,\gamma'*\gamma)$.
\end{itemize}

There is a unique $2$-functor
\begin{equation*}
\C{Cat}\To\C{Cat}_2
\end{equation*}
which is the identity on objects and $1$-morphisms. Moreover,
there is also a unique $2$-functor
\begin{equation*}
\C{Cat}\To\C{Cat}_\F
\end{equation*}
which is the identity on objects and sends a functor
$\varphi\colon\C{C}\r\C{D}$ to the identity natural transformation
$1_\varphi\colon\varphi\Rightarrow\varphi$ regarded as a morphism
$1_\varphi\colon\C{C}\r\C{D}$ in $\C{Cat}_\F$.

\begin{prop}\label{factor}
There is defined a $2$-functor $\F\colon\C{Cat}_\F\r\C{Cat}_2$
fitting into a commutative diagram
$$\xymatrix{\C{Cat}\ar[r]^\F\ar[d]&\C{Cat}\ar[d]\\\C{Cat}_\F\ar[r]_\F&\C{Cat}_2}$$
where the vertical arrows are the $2$-functors previously defined.
\end{prop}

\begin{proof}
The new $2$-functor $\F$ is defined in the following way, we use
the notation in Definition \ref{fc} and (3.A):
\begin{itemize}
\item[] \textit{on objects:} $\F\C{C}$ is the factorization
category,

\medskip

\item[] \textit{on $1$-morphisms:} the functor
$\F(\alpha)\colon\F\C{C}\r\F\C{D}$ is defined

\medskip

\begin{itemize}
\item[] \textit{on objects:} given an object $f$ in $\F\C{C}$, which is a morphism $f\colon
X\r Y$ in $\C{C}$,
$\F(\alpha)(f)=\alpha_Y\varphi(f)=\psi(f)\alpha_X$;

\smallskip

\item[] \textit{on morphisms:}
$\F(\alpha)(h,k)=(\varphi(h),\psi(k))$.
\end{itemize}

\medskip

\item[] \textit{on $2$-morphisms:}
$\F(\varepsilon,\gamma)\colon\F(\alpha)\Rightarrow\F(\beta)$ is
the natural transformation which evaluated on $f$ as above is the
morphism $\F(\varepsilon,\gamma)_f=(\varepsilon_X,\gamma_Y)$ in
$\F\C{D}$.
\end{itemize}

It is a straightforward exercise to check that this definition is
consistent and $\F$ is indeed a $2$-functor. Moreover, the diagram
in the statement commutes because all $2$-functors are the
identity on objects, $\F(1_\varphi)(f)=\varphi(f)=\F(\varphi)(f)$
and $\F(1_\varphi)(h,k)=(\varphi(h),\varphi(k))=\F(\varphi)(h,k)$.
\end{proof}

\section{Baues-Wirsching cohomology of categories}

\begin{defn}\label{cc}
Recall from \cite{csc} that a \emph{natural system} on $\C{C}$ is
a functor $D\colon \F\C{C}\r\Ab$. The \emph{Baues-Wirsching
complex} $F^*(\C{C},D)$ of a small category $\C{C}$ with
coefficients in a natural system $D$ on $\C{C}$ is a cochain
complex of abelian groups concentrated in non-negative dimensions.
In dimension $n$ this complex is given by the following product
indexed by all sequences of morphisms of length $n-1$ in $\C{C}$
\begin{equation*}
F^n(\C{C},D)=\prod_{\bullet\st{\sigma_1}\l\cdots\st{\sigma_n}\l\bullet}D(\sigma_1\cdots\sigma_n).
\end{equation*}
In this formula we assume that a sequence of length $0$ is an
object $X$ in $\C{C}$ which we also identify with the identity
morphism $1_X$. The coordinate of $c\in F^n(\C{C},D)$ in
$\bullet\st{\sigma_1}\l\cdots\st{\sigma_n}\l\bullet$ will be
denoted by $c(\sigma_1,\dots,\sigma_n)$. The differential $d$ is
defined as
\begin{eqnarray*}
  d(c)(\sigma_1,\dots,\sigma_{n+1}) &=& D(1,\sigma_1)c(\sigma_2,\dots,\sigma_{n+1}) \\
   && +\sum_{i=1}^{n}(-1)^ic(\sigma_1,\dots,\sigma_i\sigma_{i+1},\dots,\sigma_{n+1}) \\
   && +(-1)^{n+1}D(\sigma_{n+1},1)c(\sigma_1,\dots,\sigma_n).
\end{eqnarray*}
over an $n$-cochain $c$ for $n\geq 1$, and
$d(c)(\sigma)=D(1,\sigma)c(X)-D(\sigma,1)c(Y)$ for $n=0$ and
$\sigma\colon X\r Y$.

The cohomology of $\C{C}$ with coefficients in $D$ is the
cohomology of the complex $F^*(\C{C},D)$, it is denoted by
$H^*(\C{C},D)$.
\end{defn}

Baues and Wirsching noticed that $H^*$ and $F^*$ are functors in
the category $\C{Nat}$ defined as follows:
\begin{itemize}
\item[] \textit{objects:} are pairs $(\C{C},D)$ where $D$ is a natural system on
$\C{C}$,

\smallskip

\item[] \textit{morphisms:} $(\varphi,t)\colon
(\C{C},D)\r(\C{D},E)$ are pairs given by a functor
$\varphi\colon\C{D}\r\C{C}$ and a natural transformation $t\colon
D\F(\varphi)\Rightarrow E$, and composition is given by the
formula $(\psi,s)(\varphi,t)=(\varphi\psi,s(t*1_{\F(\psi)}))$.
\end{itemize}

Let $\C{Cochain}$ be the category of cochain complexes of abelian
groups and cochain homomorphisms. As a functor
$$F^*\colon\C{Nat}\To\C{Cochain}$$ is defined as follows
\begin{itemize}
\item[] \textit{on objects:} $F^*(\C{C},D)$ is the Baues-Wirsching
complex;

\smallskip

\item[] \textit{on morphisms:}
$F^n(\varphi,t)(c)(\sigma_1,\dots,\sigma_n)=t_\sigma(c(\varphi(\sigma_1),\dots,\varphi(\sigma_n)))$,
where $\sigma=\sigma_1\cdots\sigma_n$;
\end{itemize}
and $$H^n=H^nF^*\colon\C{Nat}\To\C{Ab},\;\;n\in\Z.$$

\section{The Baues-Wirsching complex as a 2-functor}

This section is the core of the paper. Its main goal is to extend
$F^*$ to a $2$-functor from an adequate $2$-category $\C{Nat}_\F$
with the same objects as $\C{Nat}$ to the following $2$-category
$\C{Cochain}_2$:
\begin{itemize}
\item[] \textit{objects:} cochain complexes of abelian groups;

\smallskip

\item[] \textit{$1$-morphisms:} cochain homomorphisms, that is
graded homomorphisms $p\colon A^*\r B^*$ of degree $0$ such that
$dp=pd$;

\smallskip

\item[] \textit{$2$-morphisms:} $[h]\colon p\Rightarrow q$ are
relative homotopy classes of homotopies between $p$ and $q$, that
is $[h]$ is represented by a degree $-1$ homomorphism $h\colon
A^*\r B^*$ such that $dh+hd=-p+q$ and $[h]=[h']$ if there exists
$r\colon A^*\r B^*$ of degree $-2$ such that $dr-rd=-h+h'$;
vertical composition of $2$-morphisms is given by
$[h'][h]=[h'+h]$, and the horizontal composition of $[h]$ and
$[h']$ in the following diagram
$$\xymatrix@C=50pt{A^*\ar@/^12pt/[r]^p_{\;}="a"\ar@/_12pt/_q[r]^{\;}="b" &B^*\ar@{=>}^{\;h}"a";"b"
\ar@/^12pt/[r]^{p'}_{\;}="c"\ar@/_12pt/[r]_{q'}^{\;}="d"
&C^*\ar@{=>}^{\;h'}"c";"d"}$$ is $[h']*[h]=[h'p+q'h]=[p'h+h'q]$;
one can use the degree $-2$ homomorphism $h'h\colon A^*\r C^*$ to
check the last equality.
\end{itemize}
Notice that morphism categories in $\C{Cochain}_2$ are in fact
groupoids. Moreover, there is a unique $2$-functor
\begin{equation}\label{2cochain}
\imath\colon\C{Cochain}\To\C{Cochain}_2
\end{equation}
which is the identity on objects and $1$-morphisms.

Let us define the $2$-category $\C{Nat}_\F$:
\begin{itemize}
\item[(5.B)\hspace{8pt}] \textit{objects:} are pairs $(\C{C},D)$ where $D$ is a natural system on
$\C{C}$;

\smallskip

\item[] \textit{$1$-morphisms:}
$(\alpha,t)\colon(\C{C},D)\r(\C{D},E)$ are pairs given by a
natural transformation $\alpha\colon\varphi\Rightarrow\psi$
between functors $\varphi,\psi\colon\C{D}\r\C{C}$, or equivalently
a morphism $\alpha\colon\C{D}\r\C{C}$ in $\C{Cat}_\F$, see (3.A),
and a natural transformation $t\colon D\F(\alpha)\Rightarrow E$,
where $\F$ is the functor defined in Proposition \ref{factor}, and
composition is defined as
$(\beta,s)(\alpha,t)=(\alpha*\beta,s(t*1_{\F(\beta)}))$;

\smallskip

\item[] \textit{$2$-morphisms:}
$(\varepsilon,\gamma)\colon(\alpha,t)\Rightarrow(\beta,s)$ are
$2$-morphisms $(\varepsilon,\gamma)\colon\alpha\Rightarrow\beta$
in $\C{Cat}_\F$ such that $t=s(1_D*\F(\varepsilon,\gamma))$, that
is the following diagram of natural transformations commutes
\setcounter{equation}{2}
\begin{equation}\label{oqc}
\xymatrix@R=15pt{D\F(\alpha)\ar@{=>}[rd]^t\ar@{=>}[dd]_{1_D*\F(\varepsilon,\gamma)}&\\&E\\D\F(\beta)\ar@{=>}[ru]_s&}
\end{equation} vertical and horizontal compositions of
$2$-morphisms in $\C{Nat}_\F$ are defined as in $\C{Cat}_\F$, that
is
$(\varepsilon',\gamma')(\varepsilon,\gamma)=(\varepsilon\varepsilon',\gamma'\gamma)$
and given a diagram in $\C{Nat}_\F$
$$\xymatrix@C=70pt{(\C{C},D)\ar@/^17pt/[r]^{(\alpha,t)}_{\;}="a"\ar@/_17pt/_{(\beta,s)}[r]^{\;}="b" &(\C{D},E)\ar@{=>}^{\;(\varepsilon,\gamma)}"a";"b"
\ar@/^17pt/[r]^{(\alpha',t')}_{\;}="c"\ar@/_17pt/[r]_{(\beta',s')}^{\;}="d"
&(\C{E},G)\ar@{=>}^{\;(\varepsilon',\gamma')}"c";"d"}$$ the
horizontal composition
$(\varepsilon',\gamma')*(\varepsilon,\gamma)=(\varepsilon*\varepsilon',\gamma*\gamma')$
in $\C{Nat}_\F$ coincides with the horizontal composition of the
following diagram in $\C{Cat}_\F$
$$\xymatrix@C=70pt{\C{C}\ar@/^15pt/@{<-}[r]^{\alpha}_{\;}="a"\ar@/_15pt/@{<-}_{\beta}[r]^{\;}="b" &\C{D}\ar@{=>}^{\;(\varepsilon,\gamma)}"a";"b"
\ar@/^15pt/@{<-}[r]^{\alpha'}_{\;}="c"\ar@/_15pt/@{<-}[r]_{\beta'}^{\;}="d"
&\C{E}\ar@{=>}^{\;(\varepsilon',\gamma')}"c";"d"}$$
\end{itemize}

It is tedious but straightforward to check that $\C{Nat}_\F$ is
indeed a well-defined $2$-category. Moreover, there is a unique
$2$-functor
\begin{equation*}
\jmath\colon\C{Nat}\To\C{Nat}_\F
\end{equation*}
which is the identity on objects and sends a morphism
$(\varphi,t)$ to $(1_\varphi,t)$. This makes sense because of the
commutativity of the diagram in Proposition \ref{factor}.

\begin{thm}\label{2F}
There is defined a $2$-functor
$F^*\colon\C{Nat}_\F\r\C{Cochain}_2$ fitting into a commutative
diagram
$$\xymatrix{\C{Nat}\ar[r]^<(.3){F^*}\ar[d]_\jmath&\C{Cochain}\ar[d]^\imath\\\C{Nat}_\F\ar[r]_<(.25){F^*}&\C{Cochain}_2}$$
\end{thm}

\begin{proof}
The new $2$-functor $F^*$ is defined as follows, we use the
notation in (3.B) and (5.B):
\begin{itemize}
\item[] \textit{on objects:} $F^*(\C{C},D)$ is the Baues-Wirsching
complex;

\bigskip

\item[] \textit{on $1$-morphisms:} $$F^*(\alpha,t)(c)(\sigma_1,\dots,\sigma_n)=
t_{\sigma}(1_D*\F(1_\varphi,\alpha))_{\sigma}c(\varphi(\sigma_1),\dots,\varphi(\sigma_n)),$$
here $\sigma=\sigma_1\cdots\sigma_n$ and the formula makes sense
because $$c(\varphi(\sigma_1),\dots,\varphi(\sigma_n))\in
D(\varphi(\sigma))=(D\F(1_\varphi))(\sigma);$$

\bigskip

\item[] \textit{on $2$-morphisms:}
$F^*(\varepsilon,\gamma)=[h_{(\varepsilon,\gamma)}]\colon
F^*(\alpha,t)\r F^*(\beta,s)$ where for an $(n+1)$-dimensional
cochain $c$ if $n>0$ $h_{(\varepsilon,\gamma)}(c)$ is defined as
\begin{equation*}
\begin{array}{l}
h_{(\varepsilon,\gamma)}(c)(\sigma_1,\dots,\sigma_n)=\\
\;\;\;\;\;\;\;\;\;\;\;s_\sigma(1_D*\F(1_\xi,\gamma\alpha))_\sigma
\sum_{i=0}^n(-1)^ic(\varphi(\sigma_1),\dots,\varphi(\sigma_i),\varepsilon_{X_i},\xi(\sigma_{i+1}),\dots,\xi(\sigma_n)),
\end{array}
\end{equation*}
where $X_i$ is the source of $\sigma_i$ and/or the target of
$\sigma_{i+1}$, notice that
$$c(\varphi(\sigma_1),\dots,\varphi(\sigma_i),\varepsilon_{X_i},\xi(\sigma_{i+1}),\dots,\xi(\sigma_n))\in(D\F(\varepsilon))(\sigma);$$
and if $n=0$
$$h_{(\varepsilon,\gamma)}(c)(X)=s_{1_X}(1_D*\F(1_\xi,\gamma\alpha))_{1_X}c(\varepsilon_X).$$
\end{itemize}

A tedious but straightforward computation shows that indeed
$$dh_{(\varepsilon,\gamma)}+h_{(\varepsilon,\gamma)}d=-F^*(\alpha,t)+F^*(\beta,s).$$
For this, essentially, one only needs to use the naturality
property of natural transformations and the commutativity of
(\ref{oqc}).

It is easy to see that $F^*$ preserves composition of
$1$-morphisms.

In order to check that $F^*$ preserves vertical composition of
$2$-morphisms we consider a diagram in $\C{Nat}_\F$
$$\xymatrix@C=80pt{(\C{C},D)\ar@/^30pt/[r]^{(\alpha,t)}_{\;}="a"\ar[r]^<(.2){(\alpha',t')}^{\;}="b"_{\;}="c"\ar@/_30pt/[r]_{(\beta,s)}^{\;}="d"&
(\C{D},E)\ar@{=>}^{(\varepsilon,\gamma)}"a";"b"
\ar@{=>}^{(\varepsilon',\gamma')}"c";"d"}$$ where
$$\xymatrix{\psi\ar@{=>}[r]^\gamma&\psi'\ar@{=>}[r]^{\gamma'}&\zeta\\
\varphi\ar@{=>}[u]_\alpha&\varphi'\ar@{=>}[u]_{\alpha'}\ar@{=>}[l]_\varepsilon&\xi\ar@{=>}[u]_\beta\ar@{=>}[l]_{\varepsilon'}}$$
is a commutative diagram of natural transformations between
functors
$$\varphi,\varphi',\psi,\psi',\xi,\zeta\colon\C{D}\r\C{C}.$$
We define a degree $-2$ homomorphism
$$r_{(\varepsilon',\gamma');(\varepsilon,\gamma)}\colon F^*(\C{C},D)\To
F^*(\C{D},E)$$ in the following way, if $c$ is an $(n+2)$-cochain
with $n>0$ then
\begin{equation*}
\begin{array}{r}
r_{(\varepsilon',\gamma');(\varepsilon,\gamma)}(c)(\sigma_1,\dots,\sigma_n)=s_\sigma(1_D*\F(1_\xi,\gamma'\gamma\alpha))_\sigma
\sum_{i=0}^n\sum_{j=i}^n(-1)^{i+j}c(\varphi(\sigma_1),\dots,\\
\varphi(\sigma_i),\varepsilon_{X_i},\varphi'(\sigma_{i+1}),\dots,\varphi'(\sigma_j),\varepsilon'_{X_j},\xi(\sigma_{j+1}),\dots,\xi(\sigma_n))
\end{array}
\end{equation*}
and for $n=0$
$$r_{(\varepsilon',\gamma');(\varepsilon,\gamma)}(c)(X)=s_{1_X}(1_D*\F(1_\xi,\gamma'\gamma\alpha))_{1_X}c(\varepsilon_X,\varepsilon'_X).$$
It is hard but straightforward to check that
$$dr_{(\varepsilon',\gamma');(\varepsilon,\gamma)}-r_{(\varepsilon',\gamma');(\varepsilon,\gamma)}d=
-h_{(\varepsilon,\gamma)}-h_{(\varepsilon',\gamma')}+h_{(\varepsilon\varepsilon',\gamma'\gamma)},$$
therefore
$F^*(\varepsilon',\gamma')F^*(\varepsilon,\gamma)=F^*(\varepsilon\varepsilon',\gamma'\gamma)$.

Let us see that $F^*$ preserves horizontal composition of
$2$-morphisms. Consider a diagram in $\C{Nat}_\F$
$$\xymatrix@C=70pt{(\C{C},D)\ar@/^17pt/[r]^{(\alpha,t)}_{\;}="a"\ar@/_17pt/_{(\beta,s)}[r]^{\;}="b" &(\C{D},E)\ar@{=>}^{\;(\varepsilon,\gamma)}"a";"b"
\ar@/^17pt/[r]^{(\alpha',t')}_{\;}="c"\ar@/_17pt/[r]_{(\beta',s')}^{\;}="d"
&(\C{E},G)\ar@{=>}^{\;(\varepsilon',\gamma')}"c";"d"}$$ Here
\begin{equation*}
\begin{array}{ccc}
\xymatrix{\psi\ar@{=>}[r]^\gamma&\zeta\\\varphi\ar@{=>}[u]^\alpha&\xi\ar@{=>}[u]_\beta\ar@{=>}[l]^\varepsilon}&
\xymatrix{\ar@{}[rd]|{\text{\normalsize and}}&\\&}&
\xymatrix{\psi'\ar@{=>}[r]^{\gamma'}&\zeta'\\\varphi'\ar@{=>}[u]^{\alpha'}&\xi'\ar@{=>}[u]_{\beta'}\ar@{=>}[l]^{\varepsilon'}}
\end{array}
\end{equation*}
are commutative diagrams of natural transformations between
functors $$\varphi,\psi,\xi,\zeta\colon\C{D}\r\C{C}\;\;\text{ and
}\;\;\varphi',\psi',\xi',\zeta'\colon\C{E}\r\C{D}.$$ We define a
degree $-2$ homomorphism
$$r'_{(\varepsilon',\gamma');(\varepsilon,\gamma)}\colon F^*(\C{C},D)\To
F^*(\C{E},G)$$ over an $(n+2)$-cochain $c$ with $n>0$ as
\begin{equation*}
\begin{array}{l}
r'_{(\varepsilon',\gamma');(\varepsilon,\gamma)}(c)(\sigma_1,\dots,\sigma_n)=\\s'_\sigma(s*1_{\F(\beta')})_\sigma(1_D*\F(1_{\xi\xi'},(\gamma*\gamma')(\alpha*\alpha')))_\sigma\sum_{i=0}^n\sum_{j=i}^n(-1)^{i+j}
c(\varphi\varphi'(\sigma_1),\dots,\\\hspace{45pt}\varphi\varphi'(\sigma_i),\varphi(\varepsilon'_{X_i}),\varphi\xi'(\sigma_{i+1}),\dots,
\varphi\xi'(\sigma_j),\varepsilon_{\xi'(X_j)},\xi\xi'(\sigma_{j+1}),\dots,\xi\xi'(\sigma_n))
\end{array}
\end{equation*}
and for $n=0$
$$r'_{(\varepsilon',\gamma');(\varepsilon,\gamma)}(c)(X)=
s'_{1_X}(s*1_{\F(\beta')})_{1_X}(1_D*\F(1_{\xi\xi'},(\gamma*\gamma')(\alpha*\alpha')))_{1_X}c(\varphi(\varepsilon'_X),\varepsilon_{\xi'(X)}).$$

After a laborious computation one can check that
$$dr'_{(\varepsilon',\gamma');(\varepsilon,\gamma)}-r'_{(\varepsilon',\gamma');(\varepsilon,\gamma)}d=
-h_{(\varepsilon',\gamma')}{F}^*(\alpha,t)-{F}^*(\beta',s')h_{(\varepsilon,\gamma)}+h_{(\varepsilon*\varepsilon',\gamma*\gamma')},$$
hence
${F}^*(\varepsilon',\gamma')*{F}^*(\varepsilon,\gamma)={F}^*(\varepsilon*\varepsilon',\gamma*\gamma')$.

The commutativity of the diagram in the statement follows easily
from the commutativity of the diagram in Proposition \ref{factor}.
\end{proof}

The set $\pi_0\C{C}$ of connected components of a small category
$\C{C}$ is formed by equivalence classes $\set{X}$ of objects in
$\C{C}$. Two objects $X, Y$ are equivalent $\set{X}=\set{Y}$ if
there exists a sequence of (non-composable) morphisms in $\C{C}$
connecting them
$$X\rightarrow\bullet\leftarrow\cdots\rightarrow\bullet\leftarrow Y.$$
This defines a product-preserving functor from small categories to
sets
$$\pi_0\colon\C{Cat}\To\C{Set}$$
with $\pi_0(\varphi)\set{X}=\set{\varphi(X)}$. Moreover, one can
obtain an ordinary category $\C{M}^0$ from a $2$-category $\C{M}$
by taking $\pi_0$ on morphism categories and also an ordinary
functor $\rho^0\colon\C{M}^0\r\C{N}^0$ from a $2$-functor
$\rho\colon\C{M}\r\C{N}$. If $\C{M}$ is a category regarded as a
$2$-category with only the trivial $2$-morphisms then
$\C{M}^0=\C{M}$.

The homotopy category of cochain complexes
$\C{Cochain}/\!\simeq\;$ coincides with $\C{Cochain}_2^0$ and the
$2$-functor $\imath\colon\C{Cochain}\r\C{Cochain}_2$ in
(\ref{2cochain}) induces the natural projection
$\imath^0\colon\C{Cochain}\r\C{Cochain}/\!\simeq\;$ onto the
quotient category.

By Theorem \ref{2F} there is a commutative diagram of functors
\begin{equation}\label{quot}
\xymatrix{\C{Nat}\ar[r]^<(.3){F^*}\ar[d]_{\jmath^0}&\C{Cochain}\ar[d]^{\imath^0}\\\C{Nat}_\F^0\ar[r]_<(.25){(F^*)^0}&\C{Cochain}/\!\simeq}
\end{equation}

\begin{prop}\label{full}
The functor $\jmath^0\colon\C{Nat}\r\C{Nat}_\F^0$ is full.
\end{prop}

\begin{proof}
Let $(\alpha,t)\colon(\C{C},D)\r(\C{D},E)$ be a morphism in
$\C{Nat}_\F$, as in (5.B). We consider the new morphism
$(1_\varphi,t(1_D*\F(1_\varphi,\alpha)))\colon(\C{C},D)\r(\C{D},E)$.
Notice that
$(1_\varphi,t(1_D*\F(1_\varphi,\alpha)))=\jmath(\varphi,t(1_D*\F(1_\varphi,\alpha)))$.
There is a $2$-morphism in $\C{Nat}_\F$
$$(1_\varphi,\alpha)\colon(1_\varphi,t(1_D*\F(1_\varphi,\alpha)))\Rightarrow(\alpha,t),$$
therefore
$$\set{(\alpha,t)}=\set{\jmath(\varphi,t(1_D*\F(1_\varphi,\alpha)))}=\jmath^0(\varphi,t(1_D*\F(1_\varphi,\alpha))).$$
\end{proof}

Proposition \ref{full} shows that $\C{Nat}_\F^0$ is a quotient
category of $\C{Nat}$, because $\jmath^0$ is the identity on
objects, and diagram (\ref{quot}) yields a factorization of the
Baues-Wirsching cohomology functors through $\C{Nat}_\F^0$.

\begin{cor}\label{elcor}
The functors $H^n\colon\C{Nat}\r\C{Ab}$ $(n\in\Z)$ factor through
the quotient category $\C{Nat}_\F^0$.
\end{cor}

\section{(Co)localization and cohomology}

In this section we prove two theorems, the localization and
colocalization theorems, which allow to easy computations in
Baues-Wirsching cohomology of small categories when the
coefficient natural system satisfy some (co)locality conditions
that will be made precise later.

Here we adopt the following definition for a localization.

\begin{defn}\label{loc}
A \emph{localization} is an adjoint pair ($\varphi$ left adjoint
to $\psi$)
\begin{equation*}
\xymatrix{\C{C}\ar@<.5ex>[r]^\varphi&\C{D}\ar@<.5ex>[l]^\psi}
\end{equation*}
such that $\varphi\psi=1_\C{D}$ and the counit is the identity
natural transformation.
\end{defn}

\begin{rem}\label{otras}
\begin{enumerate}
\item A localization can also be defined as an idempotent functor
$\xi\colon\C{C}\r\C{C}$ together with a natural transformation
$\alpha\colon 1_\C{C}\Rightarrow\xi$ such that
$1_\xi*\alpha=1_\xi=\alpha*1_\xi$. This and the former definition
are equivalent. On one hand, in this situation we can take $\C{D}$
to be the image of $\xi$, $\psi\colon\C{D}\r\C{C}$ the (faithful)
inclusion functor and $\varphi\colon\C{C}\r\C{D}$ the unique
(full) functor such that $\xi=\psi\varphi$, then we obtain a
localization in the sense of Definition \ref{loc} being $\alpha$
the unit of the adjunction. On the other hand, if we are as in the
previous definition the functor $\xi=\psi\varphi$ is idempotent
and the unit of the adjunction $\alpha\colon 1\Rightarrow\xi$
satisfies $1_\xi*\alpha=1_\xi=\alpha*1_\xi$.

\item A related concept is that of a category of fractions. Given a
set $\Sigma$ of morphisms in $\C{C}$ the \emph{category of
fractions} $\C{C}[\Sigma^{-1}]$ is a category together with a
functor $\C{C}\r\C{C}[\Sigma^{-1}]$ which is universal among all
functors which send morphisms in $\Sigma$ to isomorphisms. Let
$\C{C}_\Sigma$ be the full subcategory of $\C{C}$ formed by those
objects $X$ such that $\C{C}(f,X)$ is a bijection for any
$f\in\Sigma$ and $\psi\colon\C{C}_\Sigma\r\C{C}$ the full
inclusion. If there exists a functor
$\varphi\colon\C{C}\r\C{C}_\Sigma$ as in Definition \ref{loc} then
the composite $\C{C}_\Sigma\r\C{C}\r\C{C}[\Sigma^{-1}]$ is known
to be an equivalence of categories. On the other hand if we are in
the situation of Definition \ref{loc} and $\Sigma$ is the set of
morphisms $f$ such that $\varphi(f)$ is an isomorphism then $\psi$
induces an equivalence between $\C{D}$ and $\C{C}_\Sigma$.
\end{enumerate}
\end{rem}

\begin{defn}\label{loc2}
In the situation of Definition \ref{loc}, a natural system $D$ on
$\C{C}$ is said to be \emph{$\C{D}$-local} if given $f\colon X\r
Y$ in $\C{C}$ such that $\varphi(f)$ is an isomorphism then
$D(1_X,f)\colon D(1_X)\r D(f)$ is also an isomorphism.
\end{defn}

The following proposition gives alternative characterizations of
$\C{D}$-local natural systems.

\begin{prop}\label{loc3}
Given a natural system $D$ on $\C{C}$ the following conditions are
equivalent:
\begin{enumerate}
\item $D$ is $\C{D}$-local,

\item $D$ is naturally isomorphic to $E\F(\alpha)$
for some natural system $E$ on $\C{C}$,

\item $D$ is naturally isomorphic to $D\F(\alpha)$.
\end{enumerate}
Here $\alpha\colon 1_\C{C}\Rightarrow\psi\varphi$ is the unit
of the adjunction.
\end{prop}

\begin{proof}
Let $f\colon X\r Y$ be a morphism in $\C{C}$ such that
$\varphi(f)$ is an isomorphism, then
$E\F(\alpha)(1_X,f)=E(1_X,\psi\varphi(f))$ is an isomorphism
because $(1_X,\psi\varphi(f))$ is an isomorphism in $\F\C{C}$ with
inverse $(1_X,\psi(\varphi(f)^{-1}))$.

Suppose now that $D$ is $\C{D}$-local. The natural transformation
$$1_D*\F(1_{1_\C{C}},\alpha)\colon D=D\F(1_{1_\C{C}})\To D\F(\alpha)$$ is a natural
isomorphism since given an object $f$ in $\F\C{C}$, which is a
morphism $f\colon X\r Y$ in $\C{C}$, we have that
$(1_D*\F(1_{1_\C{C}},\alpha))_f=D(1_X,\alpha_Y)$ and
$\varphi(\alpha_Y)=1_{\varphi(Y)}$. For this last equality we use
that $\varphi$ is left-adjoint to $\psi$, $\alpha$ is the unit,
and $\varphi\psi=1$.
\end{proof}

\begin{thm}\label{locthm}
Under the conditions of Definition \ref{loc}, if $D$ is a
$\C{D}$-local natural system on $\C{C}$ then $\psi$ induces
isomorphisms $(n\in\Z)$
$$H^n(\C{C},D)\st{\simeq}\To H^n(\C{D},D\F(\psi)).$$
\end{thm}

\begin{proof}
The morphism in the statement is $H^n(\psi,1_{D\F(\psi)})$. If
$\alpha\colon1_\C{C}\Rightarrow\psi\varphi$ is the unit we have
that $\alpha*\alpha=\alpha$, see Remark \ref{otras} (1) and the
beginning of the proof of the following lemma, hence by
Proposition \ref{loc3} we can suppose without loss of generality
that $D=D\F(\alpha)$. Consider the natural transformation
$$1_D*\F(\alpha,1_{\psi\varphi})\colon
D\F(\psi)\F(\varphi)=D\F(\psi\varphi)=D\F(1_{\psi\varphi})\Rightarrow
D\F(\alpha)=D.$$ We claim that
$H^n(1_\varphi,1_D*\F(\alpha,1_{\psi\varphi}))$ is the inverse of
$$H^n(\psi,1_{D\F(\psi)})=H^n(1_\psi,1_{D\F(1_\psi)}),$$ for this
equality we use the commutativity of the diagram in Proposition
\ref{2F}.

On one hand we have the following equalities in $\C{Nat}_\F$
\begin{eqnarray*}
(1_\psi,1_{D\F(1_\psi)})(1_\varphi,1_D*\F(\alpha,1_{\psi\varphi}))&=&(1_\varphi*1_\psi,1_{D\F(1_\psi)}(1_D*\F(\alpha,1_{\psi\varphi})*1_{\F(1_\psi)}))\\
&=&(1_{\varphi\psi},1_D*\F(\alpha,1_{\psi\varphi})*\F(1_\psi,1_\psi))\\
&=&
(1_\C{D},1_D*\F(\alpha*1_\psi,1_{\psi\varphi}*1_\psi))\\
&=&(1_\C{D},1_D*\F(1_\psi,1_{\psi\varphi\psi}))\\
&=& (1_\C{D},1_D*\F(1_\psi,1_\psi))\\
&=& (1_\C{D},1_D*1_{\F(1_\psi)})\\
&=& (1_\C{D},1_{D\F(1_\psi)}).
\end{eqnarray*}
Here we use the equality $\alpha*1_\psi=1_\psi$ which can be
checked by using the equality $\varphi\psi=1_\C{D}$ and elementary
properties of adjoint functors.

On the other hand
\begin{eqnarray*}
(1_\varphi,1_D*\F(\alpha,1_{\psi\varphi}))(1_\psi,1_{D\F(1_\psi)})&=&
(1_\psi*1_\varphi,(1_D*\F(\alpha,1_{\psi\varphi}))(1_{D\F(1_\psi)}*1_{\F(1_\varphi)}))\\
&=&(1_{\psi\varphi},1_D*\F(\alpha,1_{\psi\varphi})).
\end{eqnarray*}

Now the theorem follows from Corollary \ref{elcor}, Remark
\ref{otras} (1) and the following lemma.
\end{proof}

\begin{lem}\label{ellema}
Let $\xi\colon\C{C}\r\C{C}$ and
$\alpha\colon1_\C{C}\Rightarrow\xi$ be a localization in the sense
of Remark \ref{otras} (1) and $D$ any natural system on $\C{C}$,
then there is a diagram in $\C{Nat}_\F$ as follows
$$\xymatrix@C=80pt{(\C{C},D\F(\alpha))\ar@/^30pt/[r]^{1}_{\;}="a"\ar[r]^<(.2){(\alpha,1)}^{\;}="b"_{\;}="c"\ar@/_30pt/[r]_{(1_\xi,1_D*\F(\alpha,1_\xi))}^{\;}="d"&
(\C{C},D\F(\alpha))\ar@{=>}^{(1_{1_\C{C}},\alpha)}"a";"b"
\ar@{<=}^{(\alpha,1_\xi)}"c";"d"}$$
\end{lem}

\begin{proof}
We have that
$\alpha*\alpha=(\alpha*1_\xi)(1_{1_\C{C}}*\alpha)=1_\xi\alpha=\alpha$.
By using the functor $\F$ defined in Proposition \ref{factor} we
see that $\F(\alpha)\F(1_\xi)=\F(\alpha*1_\xi)=\F(1_\xi)$,
$\F(\alpha,1_\xi)\colon\F(1_\xi)\r\F(\alpha)$, and
$\F(\alpha)\F(\alpha)=\F(\alpha*\alpha)=\F(\alpha)$, therefore the
$1$-morphisms in the diagram are well defined. Moreover, the
$2$-morphisms are also well-defined because
$1_{D\F(\alpha)}=1_D*1_{\F(\alpha)}$, $1_{\F(\alpha
)}=\F(1_{1_\C{C}},1_\xi)$ and
$$\F(1_{1_\C{C}},1_\xi)*\F(1_{1_\C{C}},\alpha)=\F(1_{1_\C{C}}*1_{1_\C{C}},1_\xi*\alpha)=\F(1_{1_\C{C}},1_\xi),$$
$$\F(1_{1_\C{C}},1_\xi)*\F(\alpha,1_\xi)=\F(1_{1_\C{C}}*\alpha,1_\xi*1_\xi)=\F(\alpha,1_\xi).$$
\end{proof}

\begin{exm}
Here we exhibit some examples of the usefulness of Theorem
\ref{locthm} as a computational tool in cohomology of categories.
We recall from \cite{csc} that any functor
$\C{C}^\mathrm{op}\times\C{C}\r\Ab$ yields a natural system on
$\C{C}$ through the natural functor
$\F\C{C}\r\C{C}^\mathrm{op}\times\C{C}$ which sends an object $f$,
which is a morphism $f\colon X\r Y$ in $\C{C}$, to the pair
$(X,Y)$.
\begin{enumerate}
\item In his study of the $2$-dimensional cohomology groups of the category $\Ab_0$ of finitely generated abelian groups
Hartl \cite{sccfgag} computed
\begin{equation}\label{moore}
H^2(\Ab_0,\ext^1_\Z(-,-\otimes\Z/2))\simeq\Z/2.
\end{equation}
The generator is the characteristic class associated to the stable
homotopy category of compact Moore spaces, see \cite{csc} 3.2. Let
$\psi\colon\C{vect}(\Z/2)\r\Ab_0$ be the inclusion of the full
subcategory of finite-dimensional $\Z/2$-vector spaces, whose
left-adjoint is $\varphi=-\otimes\Z/2$. Clearly
$\ext^1_\Z(-,-\otimes\Z/2)$ is $\C{vect}(\Z/2)$-local. Moreover,
there is a natural isomorphism
$$\ext^1_\Z(-,-\otimes\Z/2)\F(\psi)\simeq\hom_{\Z/2}$$ induced by
the change of rings spectral sequence associated to the natural
projection $\Z\twoheadrightarrow\Z/2$, see \cite{ugss} 10.2 (c),
therefore by Theorem \ref{locthm}
$$H^n(\Ab_0,\ext^1_\Z(-,-\otimes\Z/2))\simeq H^n(\C{vect}(\Z/2),\hom_{\Z/2}),\;\;\;n\in\Z.$$
Moreover, by Theorem A and Corollary 3.11 in \cite{cat} this is
the Mac Lane cohomology of $\Z/2$,
$$H^n(\C{vect}(\Z/2),\hom_{\Z/2})\simeq
H^n_{ML}(\Z/2,\Z/2)\simeq\left\{%
\begin{array}{ll}
    \Z/2, & \hbox{$n$ even;} \\
     & \\
    0, & \hbox{$n$ odd;}
\end{array}%
\right.$$ which was computed in \cite{mlcff}, therefore we recover
the isomorphism (\ref{moore}) from this complete calculation of
$H^*(\Ab_0,\ext^1_\Z(-,-\otimes\Z/2))$.

\item More generally by Theorems A and B in \cite{cat} for any prime $\ell$ there is a
spectral sequence converging to $H^{n+m}_{ML}(\Z,\Z/\ell)$ with
$$E^{nm}_2=H^n(\Ab_0,\ext^m_\Z(-,-\otimes\Z/\ell)).$$ Let
$\psi\colon\C{vect}(\Z/\ell)\r\Ab_0$ be the full inclusion of the
category of finite-dimensional $\Z/\ell$-vector spaces, whose
left-adjoint is $\varphi=-\otimes\Z/\ell$. It is obvious that
$\ext^m_\Z(-,-\otimes\Z/\ell)$ is $\C{vect}(\Z/\ell)$-local for
all
$m\geq 0$, moreover $$\ext^m_\Z(-,-\otimes\Z/\ell)\F(\psi)\simeq\left\{%
\begin{array}{ll}
    \hom_{\Z/\ell}, & \hbox{for $m=0,1$;} \\
&\\
    0, & \hbox{otherwise;} \\
\end{array}%
\right.    $$ hence for $m=0,1$
$$H^n(\Ab_0,\ext^m_\Z(-,-\otimes\Z/\ell))\simeq H^n(\C{vect}(\Z/\ell),\hom_{\Z/\ell})\simeq
H^n_{ML}(\Z/\ell,\Z/\ell)$$ and zero for $m\neq 0,1$. Therefore we
recover from the previous spectral sequence a long exact sequence
$(n\in\Z)$
$$\cdots\r H^{n-2}_{ML}(\Z/\ell,\Z/\ell)\r H^n_{ML}(\Z/\ell,\Z/\ell)\r H^n_{ML}(\Z,\Z/\ell)\r H^{n-1}_{ML}(\Z/\ell,\Z/\ell)\r\cdots$$
which was used in \cite{mlcff} to compute the Mac Lane cohomology
$H^*_{ML}(\Z,\Z/\ell)$.
\end{enumerate}
\end{exm}


\begin{defn}\label{coloc}
A \emph{colocalization} is an adjoint pair ($\varphi$ right
adjoint to $\psi$)
\begin{equation*}
\xymatrix{\C{C}\ar@<.5ex>[r]^\varphi&\C{D}\ar@<.5ex>[l]^\psi}
\end{equation*}
such that $\varphi\psi=1_\C{D}$ and the unit is the identity
natural transformation.
\end{defn}

Dually to Remark \ref{otras} we have the following observations.

\begin{rem}\label{otras2}
\begin{enumerate}
\item A colocalization can also be regarded as an idempotent functor
$\xi\colon\C{C}\r\C{C}$ together with a natural transformation
$\alpha\colon \xi\Rightarrow 1_\C{C}$ such that
$1_\xi*\alpha=1_\xi=\alpha*1_\xi$. This is equivalent to the
previous definition. Definition \ref{coloc} yields this structure
by taking $\xi=\psi\varphi$ and $\alpha$ the counit. Moreover, in
the situation of this example we obtain a colocalization as in
Definition \ref{coloc} by taking $\C{D}$ as the image of $\xi$,
$\psi\colon\C{D}\r\C{C}$ the inclusion and
$\varphi\colon\C{C}\r\C{D}$ the unique functor such that
$\xi=\psi\varphi$.

\item Given a
set $\Sigma$ of morphisms in $\C{C}$ we define ${}_\Sigma\C{C}$ to
be the full subcategory of $\C{C}$ formed by those objects $X$
such that $\C{C}(X,f)$ is a bijection for any $f\in\Sigma$. Let
$\psi\colon{}_\Sigma\C{C}\r\C{C}$ be the full inclusion. If there
exists a functor $\varphi\colon\C{C}\r{}_\Sigma\C{C}$ as in
Definition \ref{coloc} then the composite
${}_\Sigma\C{C}\r\C{C}\r\C{C}[\Sigma^{-1}]$ is an equivalence of
categories. On the other hand if we are in the situation of
Definition \ref{coloc} and $\Sigma$ is the set of morphisms $f$
such that $\varphi(f)$ is an isomorphism then $\psi$ induces an
equivalence between $\C{D}$ and ${}_\Sigma\C{C}$.
\end{enumerate}
\end{rem}

\begin{defn}\label{coloc2}
Under the conditions of Definition \ref{coloc}, a natural system
$D$ on $\C{C}$ is said to be \emph{$\C{D}$-colocal} if given
$f\colon X\r Y$ in $\C{C}$ such that $\varphi(f)$ is an
isomorphism then $D(f,1_Y)\colon D(1_Y)\r D(f)$ is also an
isomorphism.
\end{defn}

\begin{prop}\label{coloc3}
Given a natural system $D$ on $\C{C}$ the following conditions are
equivalent:
\begin{enumerate}
\item $D$ is $\C{D}$-colocal,

\item $D$ is naturally isomorphic to $E\F(\alpha)$
for some natural system $E$ on $\C{C}$,

\item $D$ is naturally isomorphic to $D\F(\alpha)$.
\end{enumerate}
Here $\alpha\colon \psi\varphi\Rightarrow 1_\C{C}$ is the counit
of the adjunction.
\end{prop}

\begin{proof}
Given a morphism $f\colon X\r Y$ in $\C{C}$ such that $\varphi(f)$
is an isomorphism, then $E\F(\alpha)(f,1_Y)=E(\psi\varphi(f),1_Y)$
is an isomorphism because $(\psi\varphi(f),1_Y)$ is an isomorphism
in $\F\C{C}$ with inverse $(\psi(\varphi(f)^{-1}),1_Y)$.

Conversely if $D$ is $\C{D}$-colocal one can readily check that
the natural transformation
$$1_D*\F(\alpha,1_{1_\C{C}})\colon D=D\F(1_{1_\C{C}})\To D\F(\alpha)$$ is a natural
isomorphism, compare the proof of Proposition \ref{loc3}.
\end{proof}

\begin{thm}\label{colocthm}
Under the conditions of Definition \ref{coloc}, if $D$ is a
$\C{D}$-colocal natural system on $\C{C}$ then $\psi$ induces
isomorphisms $(n\in\Z)$
$$H^n(\C{C},D)\st{\simeq}\To H^n(\C{D},D\F(\psi)).$$
\end{thm}

The proof of this theorem is very similar to that of Theorem
\ref{locthm}, for this reason we leave it to the reader. In this
case one has to use Proposition \ref{coloc3}, Remark \ref{otras2}
(1) and the following lemma instead of Proposition \ref{loc3},
Remark \ref{otras} (1) and Lemma \ref{ellema}.

\begin{lem}
Let $\xi\colon\C{C}\r\C{C}$ and
$\alpha\colon\xi\Rightarrow1_\C{C}$ be a colocalization in the
sense of Remark \ref{otras2} (1) and $D$ any natural system on
$\C{C}$, then there is a diagram in $\C{Nat}_\F$ as follows
$$\xymatrix@C=80pt{(\C{C},D\F(\alpha))\ar@/^30pt/[r]^{1}_{\;}="a"\ar[r]^<(.2){(\alpha,1)}^{\;}="b"_{\;}="c"\ar@/_30pt/[r]_{(1_\xi,1_D*\F(1_\xi,\alpha))}^{\;}="d"&
(\C{C},D\F(\alpha))\ar@{=>}^{(\alpha,1_{1_\C{C}})}"a";"b"
\ar@{<=}^{(1_\xi,\alpha)}"c";"d"}$$
\end{lem}

The proof of this lemma is a mere verification as in Lemma
\ref{ellema}.

\begin{exm}
Here we combine the localization and colocalization theorems
(Theorems \ref{locthm} and \ref{colocthm}) to reduce the
computation of certain cohomology groups of a triangulated
category, see \cite{mha} Chapter IV for basic facts.

More precisely, let $\C{D}$ be a small triangulated category,
whose translation functor is denoted by $\Sigma$, equipped with a
$t$-structure $(\C{D}^{\leq 0}, \C{D}^{\geq 0})$ and
$\C{A}=\C{D}^{\leq 0}\cap \C{D}^{\geq 0}$ its core. We write
$t_{\leq0}\colon\C{D}\r\C{D}^{\leq 0}$ and
$t_{\geq0}\colon\C{D}\r\C{D}^{\geq 0}$ for the right and left
adjoints of the corresponding full inclusions. Notice that the
restriction of $t_{\geq0}$ to $\C{D}^{\leq0}$, which is usually
denoted by $H^0\colon\C{D}^{\leq 0}\r\C{A}$, is also
a left-adjoint to the full inclusion. 

We consider the natural system on $\C{D}$ given by the functor
\begin{equation}\label{der}
\C{D}(t_{\leq0},\Sigma^nt_{\geq0})\colon\C{D}^\mathrm{op}\times\C{D}\To\Ab.
\end{equation}
The restriction to $\C{A}^\mathrm{op}\times\C{A}$ will be denoted
by $\ext^n_\C{D}$. Moreover, it coincides with the ordinary
extension functor $\ext^n_\C{A}$ in the abelian category $\C{A}$
when $\C{D}=\mathcal{D}^b(\C{A})$ is the bounded derived category
of $\C{A}$ with the canonical $t$-structure.

By using elementary properties of triangulated categories and
$t$-structures one can check that the natural system defined by
(\ref{der}) is $\C{D}^{\leq0}$-colocal. Moreover, its restriction
to $\C{D}^{\leq0}$, that is $\C{D}^{\leq0}(-,\Sigma^nH^0)$, is
$\C{A}$-local, therefore
$$H^*(\C{D},\C{D}(t_{\leq0},\Sigma^nt_{\geq0}))\simeq
H^*(\C{D}^{\leq0},\C{D}^{\leq0}(-,\Sigma^nH^0))\simeq
H^*(\C{A},\ext^n_\C{D}).$$

This and similar examples are relevant for the approach to
homotopy theory given by the tower of categories in \cite{cfhh}.
An explicit application can be found in \cite{propera2n}.
\end{exm}

\bibliographystyle{amsplain}
\bibliography{Fernando}
\end{document}